\documentclass{amsart}
\usepackage[dvips]{color}
\usepackage{graphicx}

\usepackage{latexsym}
\usepackage{amssymb}
\usepackage{amsmath}
\usepackage{amsthm}
\usepackage{amscd}
\theoremstyle{plain}
\newtheorem{prop}{Proposition}
\newtheorem{thm}{Theorem}
\newtheorem{lem}{Lemma}
\newtheorem{Bob}{Definition}

\newtheorem{cor}{Corollary}
\theoremstyle{remark}
\newtheorem{rem}{Remark}

\newtheorem{ques}{Question}
\newcommand{\Leb}{\ensuremath{\lambda}}
\newcommand{\measure}{\ensuremath{\bold{m}_{\mathfrak{R}}}}
\newcommand{\measb}{\ensuremath{\bold{m}_{\mathfrak{R}'}}}
\begin{document}
\title{Every ergodic transformation is disjoint from almost every IET}
\author[J.\ Chaika]{Jon Chaika}

\address{Department of Mathematics, Rice University, Houston, TX~77005, USA}

\email{jonchaika@math.uchicago.edu}
\maketitle

\begin{abstract} We show that every transformation is disjoint from almost every interval exchange transformation (IET), answering a question of Bufetov. In particular, we prove that almost every pair of IETs is disjoint. It follows that the product of almost every pair is uniquely ergodic. 
A key step in the proof is showing that any sequence of density 1 contains a rigidity sequence for almost every IET, strengthening a result of Veech. 
\end{abstract}
\begin{Bob}
Given $L=(l_1,l_2,...,l_d)$
where $l_i \geq 0$,  ${l_1+...+l_d=1}$, we obtain $d$ subintervals of $[0,1)$, ${I_1=[0,l_1) },
{I_2=[l_1,l_1+l_2)},...,{I_d=[l_1+...+l_{d-1},1)}$. Given
 a permutation $\pi$ on $\{1,2,...,d\}$, we obtain a d-\emph{Interval Exchange Transformation} (IET)  $ T_{L,\pi} \colon [0,1) \to
 [0,1)$ which exchanges the intervals $I_i$ according to $\pi$. That is, if $x \in I_j$ then $$T_{L,\pi}(x)= x - \underset{k<j}{\sum} l_k +\underset{\pi(k')<\pi(j)}{\sum} l_{k'}.$$ 
\end{Bob}
When there is no cause for confusion the subscript in denoting the IET will be omitted.
Interval exchange transformations with a fixed permutation on $d$ letters are parametrized by the standard simplex in $\mathbb{R}$, $\Delta_d=\{(l_1,...,l_d): l_i \geq 0, \sum l_i=1\}$ and $\mathring{\Delta}_d$ denotes its interior, $ \{(l_1,...,l_d): l_i > 0, \sum l_i=1\}$. In this paper, $\Leb$ denotes Lebesgue measure on the unit interval.
 The term ``almost all" refers to Lebesgue measure on the disjoint union of the simplices corresponding to the permutations that contain some IETs with dense orbits. That is, $\pi(\{1,...,k\}) \neq \{1,...,k\}$ for $k<d$ \cite[Section 3]{IET}. These permutations are called \emph{irreducible}. 
 
 Throughout this paper we assume that all measure preserving transformations are invertible transformations of Lebesgue spaces.
\begin{Bob} Two measure preserving systems $(T,X,\mu)$ and $(S,Y,\nu)$ are called disjoint (or have trivial joinings) if $\mu \times \nu$ is the only invariant measure of ${T \times S \colon X \times Y \to X \times Y}$ by $(T \times S) (x,y)=(Tx,Sy)$ with projections $\mu$ and~$\nu$.
\end{Bob}
The main result of this paper is the following theorem.
\begin{thm} \label{disjoint} Let $T \colon X \to X$ be $\mu$  ergodic. $(T,X,\mu)$ is disjoint from almost every IET.
\end{thm}
Disjointness is a way of saying that two dynamical systems are very different. It implies that they have no common factors \cite[p.\ 127 or Theorem 8.4]{glas}. For any IET $T$, and any other IET $S$, $STS^{-1}$ is an IET measurably conjugated to $T$ and therefore every IET has uncountably many IETs with which it has nontrivial joinings.
As a consequence of Theorem \ref{disjoint} we obtain a corollary.
\begin{cor} \label{ue} For any uniquely ergodic IET $T$ and almost every  IET $S$, the product $T \times S$ is uniquely ergodic. In particular, for almost every pair of IETs $(T,S)$ the product is uniquely ergodic.
\end{cor}
We prove Theorem \ref{disjoint} by the following criterion \cite[Theorem 2.1]{hp},
see also \cite[Lemma 1]{thouv} and \cite[Theorem 6.28]{glas}.
\begin{thm} (Hahn and Parry) \label{hp} If $T_1$ and $T_2$ are ergodic transformations of $(X_1, B_1, m_1)$ and $(X_2, B_2, m_2)$ respectively, and if $U_{T_1}$ and $U_{T_2}$ are spectrally singular modulo constants then $T_1$ and $T_2$ are disjoint.
\end{thm}
Recall that $U_{T_1}$ and $U_{T_2}$ are called spectrally singular modulo constants if for any functions $f \in L^2(m_1)$ and $g \in L^2(m_2)$ with integral zero, the spectral measures $\sigma_{f,T_1}$ and $\sigma_{g,T_2}$ are singular as measures. See Section 4 for a definition of spectral measures.
Spectral singularity is established by showing that for any transformation $T$ and almost every IET $S$, there exists a sequence $n_1,n_2,$...\ such that $$\underset{i \to \infty}{\lim} \int_{\mathbb{T}} z^{n_i}d \sigma_{f,S} \to \sigma_{f,S}(\mathbb{T})$$ while for any $k$ we have $$\underset{i \to \infty}{\lim}\int_{\mathbb{T}} z^{n_i+k}d \sigma_{g,T} \to 0$$ for any $f\in L^2(m_1)$ and $g\in L^2(\Leb)$ of integral 0.
To establish this result rigidity sequences are used. Given an IET $T$, a sequence $n_1,n_2,...$ is a \emph{rigidity sequence for $T$} if ${\int_0 ^1 |T^{n_i}(x)-x| d\Leb \to 0}$. This notion can be easily generalized to systems that are not IETs.
Veech proved that almost every IET has a rigidity sequence \cite[Part I, Theorem 1.3]{metric} with the following Theorem \cite[Part I, Theorem 1.4]{metric} by choosing $N_i$ corresponding to $\epsilon_i$ where $\underset{i \to \infty}{\lim} \epsilon_i=0.$
\begin{thm} \label{v rig} (Veech) For almost every interval exchange transformation $T$, with irreducible permutation, and given $\epsilon>0$ there are $N \in \mathbb{N},$ and an interval $J \subset [0,1)$ such that:
\begin{enumerate}
\item $J \cap T^n(J)=\emptyset$ for  $0<n<N$. 
\item $T$ is continuous on $T^n(J) $ for $ 0 \leq n<N$.
\item $\Leb(\underset{n=1}{\overset{N}{\cup}}T^n(J))>1-\epsilon$. 
\item $\Leb(T^N(J) \cap J)> (1-\epsilon) \Leb(J)$.
\end{enumerate}
\end{thm}
  In this paper we strengthen Veech's result that almost every IET has a rigidity sequence (see also Remark \ref{stronger} for a strengthening of Theorem \ref{v rig}).
\begin{thm} \label{rig often} Let $A$ be a sequence of natural numbers with density 1. Almost every IET has a rigidity sequence contained in $A$.
\end{thm} 

Similar classification questions have been considered in \cite{BCF}, which shows that certain pairs of 3-IETs are not isomorphic, \cite{join1} which shows that every IET is disjoint from any mixing transformation and \cite{BF} which shows that almost every IET in some permutations are disjoint from all ELF transformations. In other settings, \cite{rank 1} shows that almost every pair of rank 1 transformations is disjoint and \cite{mpd} shows that each ergodic measure preserving transformation is disjoint from a residual set of ergodic measure preserving transformations. 

 The first section provides a brief introduction to Rauzy-Veech induction and the terminology used in the second section. The results in this section are well known. The second section contains the proof of Theorem \ref{rig often}. The third section provides further consequences of the intermediate results contained in the second section. The fourth section contains the proof of Theorem \ref{disjoint}, which uses the results in the previous two sections. The final section contains consequences of Theorem \ref{disjoint} and some questions. 


\section{Rauzy-Veech Induction}

Our treatment of Rauzy-Veech induction will be the same as in \cite[Section 7]{gauss}.  We recall it here. 
 Let $T$ be a $d$-IET with permutation $\pi$. Let $\delta_+$ be the rightmost discontinuity of $T$ and $\delta_-$ be the rightmost discontinuity of $T^{-1}$. Let $\delta_{max}=\max\{\delta_+,\delta_-\}$. Consider the induced map of $T$ on $[0,\delta_{\max})$ denoted $T|_{[0,\delta_{\max})}$. If $\delta_+ \neq \delta_-$ this is a $d$-IET on a smaller interval, perhaps with a different permutation.

We can renormalize it so that it is once again a $d$-IET on $[0,1)$. That is, let $R(T)(x)= T|_{[0,\delta_{\max})}(x\delta_{max}) (\delta_{max})^{-1}$.
 This is the Rauzy-Veech induction of $T$. To be explicit the Rauzy-Veech induction map is only defined if $\delta_+ \neq \delta_-$.
 If $\delta_{max}= \delta_+$ we say the first step in Rauzy-Veech induction is $a$. In this case the permutation of $R(T)$ is given by 
\begin{equation*} \pi'(j)= \begin{cases}
 \pi (j) & \quad j \leq \pi^{-1}(d)\\ \pi(d) & \quad j=\pi^{-1}(d)+1 \\ \pi(j-1) & \quad \text{otherwise}

\end{cases}.
\end{equation*}
We keep track of what has happened under Rauzy-Veech induction by a matrix $M(T,1)$ where 
\begin{equation*} M(T,1)[ij]= \begin{cases} \delta_{i,j} & \quad j \leq \pi^{-1}(d)\\
 \delta_{i, j-1} & \quad j>\pi^{-1}(d) \text{ and } i \neq d\\
\delta_{\pi^{-1}(d),j} & \quad i=d \end{cases}. 
 \end{equation*}
 If $\delta_{max}= \delta_-$ we say the first step in Rauzy-Veech induction is $b$.
 In this case the permutation of $R(T)$ is given by 
\begin{equation*} \pi'(j)= \begin{cases}
 \pi (j) & \quad \pi(j) \leq \pi(d)\\ \pi(j)+1 & \quad \pi(d) < \pi(j) < d \\ \pi(d)+1 & \quad \pi (j)=d

\end{cases}.
\end{equation*}
 We keep track of what has happened under Rauzy-Veech induction by a matrix \begin{equation*}M(T,1)[ij]= \begin{cases} 1 & \quad i=d \text{ and }j= \pi^{-1}(d) \\ \delta_{i,j} & \quad \text{ otherwise} \end{cases}. 
\end{equation*} 
The matrices described above depend on whether the step is $a$ or $b$ and the permutation $T$ has. The following well known lemmas which are immediate calculations help motivate the definition of $M(T,1)$.
\begin{lem} \label{one step} If $R(T)=S_{L, \pi'}$  then the length vector of $T$ is a scalar multiple of $M(T,1)L$. \end{lem}
 Let $M_{\Delta}=M\mathbb{R}_d^+ \cap \mathring{\Delta}_d$. Recall $\mathring{\Delta}_d$ is the interior of the simplex in $\mathbb{R}^d$. 
\begin{lem}\label{region} An IET with lengths contained in $M(T,1)_{\Delta}$
 and permutation $\pi$ has the same first step of Rauzy-Veech induction as $T$.
 \end{lem}

We define the $n^{\text{th}}$ matrix of Rauzy-Veech induction by $$M(T,n)=M(T,n-1)M(R^{n-1}(T),1).$$ 
 It follows from Lemma \ref{region} that for an IET with length vector in $M(T,n)_{\Delta}$ and permutation $\pi$ the first $n$ steps of Rauzy-Veech induction agree with $T$. If $M$ is any matrix, $C_i(M)$ denotes the $i^{th}$ column and $C_{max}(M)$ denotes the column with the largest sum of entries. Let $|C_i(M)|$ denote the sum of the entries in the $i^{th}$ column. Versions of the following lemma are well known and we provide a proof for completeness.
\begin{lem} \label{proscribed rv} If $M(R^n(T),k)$ is a positive matrix and $L= \frac {C_i(M(T,n+k))}{|C_i(M(T,n+k))|}$ then $S_{L, \pi}$ 
 agrees with $T$ through the first $n$ steps of Rauzy-Veech induction. 
\end{lem} 
\begin{proof} By Lemma \ref{one step} the length vector for $R^m(S_{L, \pi})$ is $ \frac{C_i(M(R^m(T),n+k-m))}{|C_i(M(R^m(T),n+k-m))|}$ for any $m$ where $R^m(S_{L,\pi})$ is defined. 
  By our assumption on the positivity of $M(R^n(T),k)$ the vector $\frac{C_i(M(R^n(T),k))}{|C_i(M(R^n(T),k))|}$ is contained in $\mathring{\Delta}_d$. The lemma follows by  Lemma \ref{region} and induction. 
\end{proof}
The next definition does not appear in \cite{gauss} but is important for the next section.
\begin{Bob}
 A matrix $M$ is called $\nu$ \emph{balanced} if $\frac 1 {\nu} <\frac {|C_i(M)|}{|C_j(M)|}<\nu$ for all $i$ and $j$.
\end{Bob}
 Notice that if $M$ is $\nu $~ balanced then $|C_i(M)|>\frac {|C_{max}(M)|}{\nu}$.

\section{Proof of Theorem \ref{rig often}}
Theorem \ref{rig often} follows from the following proposition.
\begin{prop} \label{good} Let $A \subset \mathbb{N}$ be a sequence of density 1. For every $\epsilon>0$ and almost every IET $S$, there exists $n_{\epsilon} \in A$ such that $\int_0 ^1 |S^{n_{\epsilon}}(x)-x|d\lambda< \epsilon$.
\end{prop}
This proposition implies Theorem \ref{rig often} because the countable intersection of sets of full measure has full measure.

 Motivated by this proposition if ${\int_0 ^1 |T^{n}(x)-x| d\Leb < \epsilon}$ we say $n$ \emph{is an }$\epsilon$ \emph{rigidity time for} $T$. 

 Throughout this section we will assume that the IETs are in a fixed Rauzy class $\mathfrak{R}$, which contains  $d$-IETs with some irreducible permutations. Let $r$ denote the number of different permutations IETs in $\mathfrak{R}$ may have. Let $\measure$ denote Lebesgue measure on $\mathfrak{R}$ (the disjoint union of $r$ simplices in $\mathbb{R}^d$).

Proposition \ref{good} will be proved by showing that there is a particular reason for $\epsilon$~rigidity (called acceptable $\epsilon$ rigidity) that occurs often in many $P_i:=[2^i,2^{i+1}]$ (Proposition \ref{rig}) but rarely occurs for any fixed $n$ (Lemma \ref{estimate}). For every IET $S$ satisfying the Keane condition, and every $i$ there exists some $n$ such that $|C_{max}(M(S,n))| \in P_i$. In general there can be more than one such $n$.



For each of the permutations $\pi_1,...,\pi_r$ that an IET in $\mathfrak{R}$ may have, fix a finite sequence of Rauzy-Veech induction steps  $\omega_i$, which gives a positive matrix. That is each letter of $\omega_i$ will be one of the two types of Rauzy-Veech steps ($a$ or $b$) and the product of the sequence of the associated matrices starting from permutation $\pi_i$ provides a positive Rauzy-Veech matrix. Let $M(\omega_i)$ denote this matrix. Let $|\omega_i|$ denote the number of steps in $\omega_i$. 
 Let $p_i=\measure(M(\omega_i)_{\Delta})$.

\begin{Bob} We say a pair $(M,C_{max}(M))$ is \emph{acceptable} if $M=M(T,n)$, $R^{n-|\omega_i|}(T)$ has permutation $\pi_i$ and $M(R^{n-|\omega_i|}(T),
 |\omega_i|)=M(\omega_i)$ for some $1\leq i\leq r$. 

If $(M, C_{max}(M))$ is an acceptable pair then $M$ is called an \emph{acceptable matrix}.
\end{Bob}
Informally, if $M=M(T,n)$ then the pair $(M,C_{max}(M))$ is acceptable if the last steps in Rauzy-Veech induction for an IET with length vector in $M_{\Delta}$ agrees with some $\omega_i$ and the permutation of $R^{n-|\omega_i|}(T)$ is $\pi_i$.
\begin{rem}\label{how to accept} In the remainder of this section we will use the fact that if $R^n(T_{L,\pi})$ has permutation $\pi_i$ then for any IET $S$ with length vector in $(M(T_{L,\pi},n)M(\omega_i))_{\Delta}$ and permutation $\pi$ the pair $(M(S,n+|\omega_i|), C_{max}(M(S,n+|\omega_i|)))$ is acceptable.
\end{rem}

\begin{lem} \label{proportional} There exists $\nu$ such that any acceptable matrix is $\nu$ balanced. \end{lem} 
\begin{proof} Let $M_1$ be a positive matrix. Observe that if $M_2$ is a matrix with nonnegative entries then $M_2M_1$ is at worst $\underset{i,j,k}{\max}\,\frac{M_1[i,j]}{M_1[i,k]}$ balanced. Since there are only finitely many $M(\omega_i)$ and they are all positive the lemma follows. In particular, we can chose $\nu= \underset{t}{\max}\,\underset{i,j,k}{\max}\,\frac{M(\omega_t)[i,j]}{M(\omega_t)[i,k]}$.
\end{proof}

\begin{lem} \label{r acceptable}For any $d$-column $C$, $|\{M \colon (M,C) \text{ is an acceptable pair } \}|\leq r^2$.
\end{lem}
That is, any $d$-column can appear in at most $r^2$ different acceptable pairs in a given Rauzy class even if the permutations are allowed to vary in the Rauzy class. 
\begin{proof} Assume $C$ belongs to two different acceptable pairs $(M(T,n),C)$, and $(M(S,n'),C)$ where both $T$ and $S$ have the same permutation $\pi_i$ (this is an additional assumption). The acceptable sequence of steps $\omega_j$ for $T$ and $\omega_{j'}$ for $S$ are different. This is because if $\omega_j=\omega_{j'}$ then the last $|\omega_j|$ steps of Rauzy-Veech induction are the same. However, since $C=C_{max}(M(T,n))=C_{max}(M(S,n'))$ and $S$ and $T$ have the same starting permutation,  Lemma \ref{proscribed rv} implies that all but the last $|\omega_j|$ steps of Rauzy-Veech induction are the same and therefore $M(T,n)=M(S,n')$. 
There can only be $r$ such pairs (with permutation $\pi_i$) 
 because there are $r$ choices of $\omega_j$. There are $r$ choices of $\pi_i$ so the lemma follows.
\end{proof}
\begin{prop}\label{rig reas} For \measure -almost every IET $S$, the set of natural numbers
  \begin{multline*}
  \lefteqn{\{i: \text{ for some } n, \,  |C_{max}(M(S,n))| \in P_i \text{ and }} \\ (M(S,n),C_{max}(M(S,n))) \text{ is an acceptable pair } \}
\end{multline*}
has positive lower density.
\end{prop}

The following two lemmas are used in the proof of Proposition \ref{rig reas}.
\begin{lem}\label{bal often}  For \measure -almost every IET $S$, and all sufficiently large $\nu_0$, the set of natural numbers
$$G(S):=\{i: \text{for some } n, \, |C_{max}(M(S,n))| \in P_i  \text{ and } M(S,n) \text{ is } \nu_0 \text{ balanced}\}$$ 
has positive lower density.
\end{lem}
\begin{rem} It is not claimed that a positive lower density of the Rauzy-Veech induction matrices are balanced.
\end{rem}
To prove this we use an independence type result for Rauzy-Veech induction that we provide a slight reformulation of \cite[Corollary 1.7]{ker}.
\begin{prop} (Kerckhoff) \label{likely bal} Let $\mathfrak{R}$ be one of the Rauzy classes of permutations of $d$-IETs. There exist $p>0,K>1$ and $\nu_0>1$ depending only on $\mathfrak{R}$ such that for any matrix of Rauzy-Veech induction $M'=M(S,n)$ we have
\begin{multline*}
\lefteqn{\measure(\{T: \pi(T)=\pi(S), T \in M'_{\Delta} \exists m>n \text{ such that } M(T,m) \text{ is }}\\ \nu_0 \text{-balanced and }  |C_{max} (M(T,m))| <K^d|C_{max}(M')|\})> p\measure ( M'_{\Delta})\end{multline*}
\end{prop}
This proposition is useful because the constants are independent of $M'$.
\begin{proof}[Proof of Lemma \ref{bal often}]  Consider the independent $\mu$ distributed random variables $F_1,F_2,...$ where $\mu$ takes value 1 with probability $p$ and 0 with probability $1-p$ and $F_i:\Omega \to \{0,1\}$. Recall that one puts a probability measure $\mu^{\mathbb{N}}$ on $\Omega$ such that for any $k\leq n$ and $a_1,...,a_n \in \{0,1\}$ where $k$ of the $a_i$ are 1 we have $$\mu^{\mathbb{N}}(\{t \in \Omega: F_i=a_i \text{ for all } i \leq n\})=p^k(1-p)^{n-k}.$$
 By the Proposition \ref{likely bal}, given $G(S) \cap [0,N]$ the conditional probability that $N + i  \in G(S)$ for some $0<i \leq \lceil d \log_2(K) \rceil $ is at least $p$. Thus by induction on $k$, for any natural numbers $n_1,n_2,...,n_k$
\begin{multline*}\lefteqn{\measure(\{S: [n_i\lceil d \log_2(K) \rceil, (n_i+1)]\lceil d \log_2(K) \rceil] \cap G(S) \neq \emptyset \, \forall i\leq k \})} \\ \geq \mu^{\mathbb{N}}(\{t: F_{n_i}(t)=1 \, \forall i \leq k\}).
\end{multline*} Briefly, assume that we are given $n_{k+1}>n_k$ and consider all $G(S) \cap [0,n_{k+1}\lceil d\log_2(K)\rceil]$ such that $[n_i,n_i + \lceil d\log_2(K)\rceil] \cap  \in G(S) \neq \emptyset$ for each $i\leq k$. By our inductive hypothesis the measure of such $S$ is at least $p^k$.  By our previous remark at least $p$ of these $S$ have $[n_{k+1},n_{k+1} + \lceil d\log_2(K)\rceil] \cap G(S) \neq \emptyset$. 
Since, by the strong law of large numbers, for $\mu^ {\mathbb{N}}$-almost every $t$ we have $\underset{n \to \infty} \lim \frac{\underset{i=1}{\overset{n}{\sum}}F_i(t)}{n}=p$
we have that for \measure  -almost every $S$, $G(S)$ has lower density at least $\frac{p}{ \lceil d \log_2(K) \rceil}$. 
\end{proof}
\begin{lem} \label{bal} (Kerckhoff) If $M$ is $\nu_0$ balanced and $W \subset \Delta_d$ is a measurable set, then $$\frac{\measure(W)}{\measure(\Delta_d)}< \frac{\measure(MW)}{\measure(M\Delta_{d})}(\nu_0)^{-d}.$$
\end{lem}
This is \cite[Corollary 1.2]{ker}. See \cite[Section 5]{ietv} for details.
\begin{proof}[Proof of Proposition \ref{rig reas}] 
By Lemma \ref{bal} if $M(T,n)$ is $\nu_0$ balanced and $R^n(T)$ has permutation $\pi_i$ then $\frac{\measure(M(T,n)M(\omega_i)\Delta_d)}{\measure(M(T,n)\Delta_{d})} \geq \nu_0^{-d} p_i$. In words: given that $M(T,n)$ is $\nu_0$ balanced and that $R^n(T)$ has permutation $\pi_i$, the conditional probability that $(M(T,n + |\omega_i|), C_{max}(M(T,n + |\omega_i|)))$ is an acceptable pair is at least $\nu_0^{-d} p_i$.  Considering each $\pi_i$, the proposition follows analogously to Lemma \ref{bal often}.
\end{proof}
\begin{Bob} Let $S$ be an IET.
If $(M(S,n), C_{max}(M(S,n)))$ is acceptable and $m=|C_{max}(M(S,n))|$ is an $\epsilon$ rigidity time for $S$ then $m$ is called an \emph{acceptable} $\epsilon$ \emph{rigidity time for} $S$.
\end{Bob}
\begin{prop} \label{rig} For every $\epsilon>0$, \measure -almost every IET $S$, the set of natural numbers 
 $$ G_{\epsilon}(S):=\{i:  P_i \text{ contains an acceptable } \epsilon \text{ rigidity time for } S \}$$
has positive lower density.
\end{prop}
\begin{proof} Consider an IET $S_{L, \pi}=S$ such that $(M(S,n), C_k(M(S,n)))$ is an acceptable pair (in particular, $C_k(M(S,n))=C_{max}(M(S,n))$).
  For ease of notation let $M'=M(S,n)$. 
Let $W_{k,\epsilon}=\{(l_1,l_2,...,l_d): \l_i>0 \, \forall i, \, l_k>1-\frac{\epsilon}{3}, \sum l_i=1\}$. 
If $L \in W_{k,\epsilon}$ then $T_{\frac{M' L}{|M' L|}, \pi}$ has an $\epsilon$ rigidity time of $|C_k(M')|$.
 This is the reason for rigidity used to prove Theorem 1.3 and 1.4 \cite[pages 1337-1338]{metric}. If  $M'$ is acceptable then Lemma \ref{proportional} states that $M'$ is $\nu$ balanced. It then follows
by Lemma \ref{bal} 
that the proportion of $M' _{\Delta}$ which has $|C_k(M')|$ as an $\epsilon$ rigidity time is at least $\nu^{-d}\measure (W_{k,\epsilon})$. Thus if $i_1< i_2<... \in G(S)$ then the probability that  $i_f \in G_{\epsilon}(S)$ is at least $\nu^{-d}\measure(W_{k, \epsilon})$ regardless of which $i_k \in G_{\epsilon}(S)$ for $k<f$.
The proposition follows analogously to Lemma \ref{bal often}. 
\end{proof}

Before proving Proposition \ref{good} we provide the following lemmas.
\begin{lem} There exists $b \in \mathbb{R}$ such that for any $n \in \mathbb{N}$, $$|\{M \colon M \text{ is acceptable and } |C_{max}(M)|=n \}|\leq  b n^{d-1}.$$
\end{lem}
\begin{rem} The constant $b$ depends only on our Rauzy class $\mathfrak{R}$. It is not claimed that for every $n \in \mathbb{N}$ there exists an acceptable matrix $M$ with $|C_{max}(M)|=n$.
\end{rem}
\begin{proof} By Lemma \ref{r acceptable} each column $C$ can be  $|C_{max}(M)|$ for at most $r^2$ different acceptable matrices $M$. By induction on $d$, $O(n^{d-1})$ different $d$-columns with non-negative integer entries have the sum of their entries equal to $n$.
\end{proof}
\begin{lem} \label{size} (Veech) If $M$ is a matrix given by Rauzy-Veech induction, then $$\measure(M_{\Delta})=c_{\mathfrak{R}} \underset{i=1}{\overset{d}{\Pi}} |C_i(M)|^{-1}.$$
\end{lem}
This is \cite[equation 5.5]{ietv}. An immediate consequence of it is that any $\nu$ balanced Rauzy-Veech matrix $M$ has $\measure(M_{\Delta}) \leq c_{\mathfrak{R}} \nu^{d-1} |C_{max}(M)|^{-d}$. The previous two lemmas give the following result.
\begin{lem} \label{estimate} The \measure -measure of IETs that have acceptable pairs with the same $|C_{max}|$ is at most $O(|C_{max}|^{-1})$.
\end{lem}
\begin{proof}[Proof of Proposition \ref{good}] By Lemma \ref{estimate} and the fact that $A$ has density 1, $$\underset{i \to \infty}{\lim}\measure(\{T \colon \exists n \text{ with } M(T,n) \text{ acceptable and } |C_{max}(M(T,n))| \in P_i \backslash A\})=0.$$  
 Therefore, Proposition \ref{rig} implies that for any $\epsilon>0$, almost every IET has an acceptable $\epsilon$ rigidity time in $A$. In fact, almost every IET has an $\epsilon$ rigidity time in $P_i \cap A$ for a positive upper density set of $i$. 
\end{proof}
\begin{rem} \label{stronger} To be explicit, Proposition \ref{rig} shows that for any sequence A with density 1, and any $\epsilon>0$, for almost every IET the integer $N$ in Veech's Theorem \ref{v rig} can be chosen from $A$.
\end{rem}

\section{Consequences of Section 2}
In this section we glean some consequences of the proofs in the previous section. One of these (Corollary \ref{eigen}) follows from \cite[Theorem A]{AF} and is used in the proof of Theorem \ref{disjoint}. It is proven independently of \cite[Theorem A]{AF} in this section.
\begin{cor} Let $A$ be a sequence of natural numbers with density 1. A residual set of IETs has a rigidity sequence contained in $A$.
\end{cor}
\begin{proof}
Take the interior of the set $W_{k, \epsilon}$ considered in the proof of Proposition \ref{rig}. In this way one obtains that the set of IETs with an $\epsilon$ rigidity time in $A$ contains an open set of full measure (therefore dense). Intersecting over $\epsilon$ shows that a residual set of IETs has a rigidity sequence in any sequence of density 1.
\end{proof}

 The number of columns that can appear in Rauzy-Veech matrices grows at least like $u_{\mathfrak{R}}R^{d}$ where the constant $u_{\mathfrak{R}}$ depends on $\mathfrak{R}$ and $R$ is the norm of the largest column of the matrix. Briefly, in order to collect a positive measure of IETs having admissible matrices $M$, with $|C_{max}(M)| \in P_k$, Lemma \ref{size} implies that there needs of be more than $u_{\mathfrak{R}}(2^k)^d$ admissible matrices with $|C_{max}| \in P_k$. This provides a partial answer to the first question in \cite[Part II, Questions 10.7]{metric} which asks what one can say about for the growth of so called primitive IETs (IETs with rational lengths that are as close to being minimal as possible). If $N(R,\pi)$ are the number of primitive IETs with permutation $\pi$ on $d$ letters and period less than $R$ it asks what one can say about $R^{-d}N(R,\pi)$.

The next result provides a slight improvement of Theorem \ref{rig often} and uses the following definition. 
\begin{Bob} Let $S$ be an IET. We say $m$ is an \emph{expected} $\epsilon$ \emph{rigidity time} for $S$ if there exists an $n$ such that that the following two conditions are met.
\begin{enumerate}
\item $(M(S,n),C_{max}(M(S,n)))$ is acceptable and $m=|C_{max}(M(S,n))|$.
\item $C_{max}(M(S,n))=C_k(M(S,n))$ and $R^n(S)$ lies in the set $W_{k,\epsilon}$ defined in the proof of Proposition \emph{\ref{rig}}.
\end{enumerate}
\end{Bob}
Every expected $\epsilon$ rigidity time is an acceptable $\epsilon$ rigidity time. 
\begin{cor}\label{better} For every $\epsilon>0$ and Rauzy class $\mathfrak{R}$ there is a constant $a_{\mathfrak{R}}(\epsilon)<1$ such that any sequence of natural numbers $A$ with density at least $a_{\mathfrak{R}}(\epsilon)$ has a rigidity sequence for all but a \measure -measure $\epsilon$ set of IETs.
\end{cor}
\begin{proof}
First note that the set of IETs having a rigidity sequence contained in $A$ is measurable.
Let $e_{\mathfrak{R}}(\epsilon)$ denote $\measure (W_{k,\epsilon})$. Let $M=M(T_{L,\pi},n)$ be an acceptable matrix. 
By the bound on distortion in Lemma \ref{bal}, 
 the conditional probability of an IET in  $M_{\Delta}$ and permutation $\pi$ having an expected $\epsilon$ rigidity time $|C_{max}(M)|$ is proportional to $e_{\mathfrak{R}}(\epsilon)$. This uses Lemma \ref{proportional} which states that if $M$ is an acceptable matrix then $M$ is $\nu$ balanced. An analogous argument to Lemma \ref{bal often} shows that there exists $c_1>0$ such that the set
$$\{i \colon \exists \, m \in P_i \text{ which is an expected } \epsilon \text{ rigidity time for }T\}$$ has lower density at least $c_1e_{\mathfrak{R}}(\epsilon)$ for almost every $T$. 
Because $(M,C_{max}(M))$ is acceptable  Lemma \ref{estimate}  establishes that there exists $c_2>0$ (where $c_2$ is the constant from the $O(n^{-1})$) 
such that $$\measure (\{T \colon n \text{ is an expected } \epsilon \text{ rigidity time for }T\})< c_2\nu^de_{\mathfrak{R}}(\epsilon)n^{-1}$$ for all $n$. Thus, for any $\epsilon>0$ and $\delta>0$ a set of natural numbers with density $1-\delta$ contains an $\epsilon$ expected rigidity time for all but a set of IETs of measure $2 \delta \frac {c_2}{c_1}\nu^d$ and the corollary follows.
\end{proof}
\begin{rem}  Recall that $\nu$ depends on the choices of $\omega_i$ that define acceptable pairs. The constant $c_1$ depends on $\nu$.
\end{rem}

Corollary \ref{better} gives two further corollaries.
\begin{cor} Almost every IET has a rigidity sequence which is not a rigidity sequence for $\measb$-almost every IET and every $\mathfrak{R}'$.
\end{cor}
\begin{proof} It suffices to show that for any $\delta>0$ and Rauzy class $\mathfrak{R}'$ all but a set of \measure -measure $\delta$  IETs have a rigidity sequence that is not a rigidity sequence for \measb -almost every IET. 
Given $\epsilon_1,\epsilon_2 >0$ and a Rauzy class, $\mathfrak{R}'$ consider the set 
\begin{multline*}\lefteqn{A_{\mathfrak{R}'}(\epsilon_1,\epsilon_2)= \{n \colon n \text{ is an } \epsilon_1 \text{ rigidity time for a set of IETs of }}\\ \measb \text{-measure at least } \epsilon_2\}.\end{multline*}
 If $\epsilon_2>0$ and $\mathfrak{R}'$ are fixed then the density of this set goes to zero with $\epsilon_1$. To see this, observe that if $n_1$ and $n_2$ are $\epsilon$ rigidity times for $T$ then $n_1-n_2$ is a $2\epsilon$ rigidity time for $T$.
 It follows that if $\epsilon< \frac 1 2 \underset{0<n \leq M}{\min}\int |T^nx-x|d\Leb$ then $\{r+1,r+2,...,r+M\}$ can contain at most one $\epsilon$ rigidity time for $T$.
Choose $\epsilon_1(k)$ so that the (upper) density of $A_{\mathfrak{R}'}(\epsilon_1(k),\frac 1 k)$ is less than $1-a_\mathfrak{R}(\delta)$.
 By Corollary \ref{better}, all but a  \measure -measure $\delta$ set of IETs have a rigidity sequence in the complement of  $A_{\mathfrak{R}'}(\epsilon_1(k),\frac 1 k)$ (which can be shared by a set of IETs with \measb -measure at most $\frac 1 k$). Consider the  countable intersection over $k$ of these sets of \measure -measure at most $1-\delta$, which also has measure at most $1-\delta$ because the sets are nested. For each IET $T$ in this set let $n_i$ be a $\frac 1 i$ rigidity time for $T$ lying in the complement of $A_{\mathfrak{R}'}(\epsilon_1(i),\frac 1 i)$. Therefore, $n_1,n_2,...$ is a rigidity sequence for $T$ that is not a rigidity sequence for \measb -almost every IET.
\end{proof}

\begin{cor} \label{eigen} For every $\alpha \notin \mathbb{Z}$, almost every IET does not have $e^{2 \pi i \alpha}$ as an eigenvalue.
\end{cor}
We will prove this corollary independently of \cite[Theorem A]{AF}, from which it immediately follows.
 \begin{thm} (Avila and Forni) If $\pi$ is an irreducible permutation that is not a rotation, then almost every IET with permutation $\pi$ is weak mixing.
\end{thm}
The proof is split into the case of rational $\alpha$ and the case of irrational $\alpha$. If $T$ has $e^{2 \pi i \alpha}$ as an eigenvalue for some rational $\alpha \notin \mathbb{Z}$ then it is not totally ergodic. This is not the case for almost every IET \cite[Part I, Theorem 1.7]{metric}.
\begin{thm} (Veech) Almost every IET is totally ergodic.
\end{thm}
It suffices to consider irrational $\alpha$ and show that for any $\delta>0$ and $\mathfrak{R}$, the set of IETs having $e^{2 \pi i \alpha}$ as an eigenvalue has \measure -outer measure less than $\delta$. If $e^{2 \pi i \alpha}$ is an eigenvalue for $T$ then rotation by $\alpha$ is a factor of $T$. However, rigidity sequences of a transformation are also rigidity sequences for the factor. 
For every irrational $\alpha$ and $e>0$ there is a sequence of density $1-e$ that contains no rigidity sequence for rotation by $\alpha$. 
To see this, observe that if $n_1$ and $n_2$ are $\epsilon$ rigidity times for $T$ then $n_1-n_2$ is a $2\epsilon$ rigidity time for $T$.
 It follows that if $\epsilon< \frac 1 2 \underset{0<n \leq M}{\min}\int |T^nx-x|d\Leb$ then $\{k+1,k+2,...,k+M\}$ can contain at most one $\epsilon$ rigidity time for $T$. Choose $e<1-a_{\mathfrak{R}}(\delta)$ and pick a sequence of density $1-e$ containing no rigidity sequence for rotation by $\alpha$. The IETs having a rigidity sequence in this sequence have \measure -measure at least $1-\delta$ and Corollary \ref{eigen} follows.


\begin{rem}
Every sequence of density 1 contains a rigidity sequence for rotation \mbox{by $\alpha$.}
\end{rem}

\section{Proof of Theorem \ref{disjoint}}
Given a $\mu$ measure preserving dynamical system $T$, let $U_T$ be the unitary operator on $L^2(\mu)$ given by $U_T(f)=f \circ T$. Let $L^2_0$ denote the set of $L^2$ functions orthogonal to constant functions. If $f \in L^2$ let $\sigma_{f,T}$ be the spectral measure for $f$ and $U_T$, that is the unique measure on $\mathbb{T}$ such that
$$ \int_{\mathbb{T}} z^n d \sigma_{f,T}=<f, U_T^n f> \, \text{ for all } n.$$

Fix $T \colon  [0,1) \to [0,1)$, a $\mu$ ergodic transformation.
 By Theorem \ref{hp}, establishing that for any $S$ in a full measure set of IETs $\sigma_{f,T}$ is singular with respect to $\sigma_{g,S}$ for any $f \in L_0^2(\mu)$ and $g \in L_0^2(\Leb)$  establishes Theorem \ref{disjoint}.
 Let $H_{pp}$ be the closure of the subspace of $L_0^2(\mu)$ spanned by non-constant eigenfunctions of $U_T$ (where the spectral measures are atomic) and $H_c$ be its orthogonal complement (where the spectral measures are continuous).
\begin{lem}\label{point spec} If $f \in H_{pp}$ then for almost every IET $S$, $\sigma_{f,T}$ is singular with respect to $\sigma_{g,S}$ for any $g \in L_0^2(\Leb)$.
\end{lem}
\begin{proof} Let $f \in H_{pp}$. The atomic measure $\sigma_{f,T}$ is supported on the $e^{2 \pi i \alpha}$ that are eigenvalues of $U_T$. If $\sigma_{f,T}$ is nonsingular with respect to $\sigma_{g,S}$ then $U_T$ and $U_S$ share an eigenvalue (other than the simple eigenvalue 1 corresponding to constant functions). The set of eigenvalues of $U_T$ is countable because $H_{pp}$ has a countable basis of eigenfunctions.  The lemma follows from the fact that the set of IETs having a particular eigenvalue has measure zero (Corollary \ref{eigen}) and the countable union of measure zero sets has measure zero.
\end{proof}
\begin{lem} \label{cont spec} If $f \in H_{c}$ then for almost every IET $S$, $\sigma_{f,T}$ is singular with respect to $\sigma_{g,S}$ for any $g \in L_0^2(\Leb)$.
\end{lem}
To prove this lemma we use Wiener's Lemma (see e.g. \cite[Lemma 4.10.2]{bs}) and its immediate corollary.
\begin{lem} (Wiener)
 For a finite measure $\mu$ on $\mathbb{T}$ set $\hat{\mu}(k)= \int_{\mathbb{T}} z^k d\mu(z)$. ${\underset{n \to \infty}{\lim} n^{-1} \underset{k=0}{\overset{n-1}{\sum}}| \hat{\mu}(k)|^2=0}$ iff $\mu$ is continuous. \end{lem}
\begin{cor} \label{wiener} For a finite continuous measure $\mu$ on $\mathbb{T}$ there exists a density 1 sequence $A$, such that $\underset{k \in A}{\lim} \, \hat{\mu}(k)=0$.
\end{cor}
\begin{proof}[Proof of Lemma \ref{cont spec}] Decompose $H_c$ into the direct sum of mutually orthogonal $H_{f_i}$, where each $H_{f_i}$ is the cyclic subspace generated by $f_i$ under $U_T$ (and $U_{T}^{-1}=U_T^*$). By Corollary~ \ref{wiener}, for each $i$ there exists a density 1 set of natural numbers $B_i$ such that $\underset{n \in B_i}{\lim}\int_{\mathbb{T}} z^{n} d\sigma_{f_i,T}=0$. Choose $N_j$ increasing such that for each $j$ we have $\underset{n>N_j}{\inf} \frac {|B_i \cap [0,n]|}{n}> 1-2^{-j}$. Let $A_i:= \underset{j=1}{\overset{\infty}{\cup}}\left([N_j, N_{j+1}] \cap  \underset{k=-j}{\overset{j}{\cap}} B_i+k\right)$. By construction, $(A_i-k) \backslash B_i$ is a finite set for any $k \in \mathbb{Z}$. Therefore, $\underset{n \in A_i}{\lim}\int_{\mathbb{T}} z^{n+k} d\sigma_{f_i,T}=0$ for any $k \in \mathbb{Z}$.
 Thus, for any $h \in H_{f_i}$ it follows that $\underset{n \in A_i}{\lim}\int_{\mathbb{T}} z^{k+n} d\sigma_{h,T}=0$ for any $k$. This follows from the fact that $\sigma_{h,T} \ll \sigma_{f_i,T}$, the span of $z^k$ is dense in $L_2$ and $|\int_{\mathbb{T}} z^{r} d\mu|\leq \mu(\mathbb{T})$.
 Since there are only a countable number of $H_{f_i}$, there exists a density 1 sequence $A$ such that for any $i$ and $h \in H_{f_i}$ we have that $\underset{n \in A}{\lim}\int_{\mathbb{T}} z^{k+n} d\sigma_{h,T}=0$ for any $k$.
The construction of $A$ is similar to the construction of the $A_i$. That is, pick $N_j$ such that $\underset{n>N_j}{\inf} \frac{|A_i \cap [0,n]|}{n}> 1-2^{-j}$ for any $i<j$. Let $$A= \underset{j=1}{\overset{\infty}{\cup}} [N_j,N_{j+1}]\cap A_1 \cap...\cap A_j.$$
 It follows that for any $h \in H_c$, $\underset{n \in A}{\lim}\int_{\mathbb{T}} z^{k+n} d\sigma_{h,T}=0$ for any $k$. This uses the fact that if $g_1$ and $g_2$ lie in orthogonal cyclic subspaces then $\sigma_{g_1+g_2,T}$ is $\sigma_{g_1,T}+\sigma_{g_2,T}$. 

Let $S$ be any IET with a rigidity sequence contained in $A$, which almost every IET has by Theorem \ref{rig often}. Notice that since $n_1,n_2,...$ is a rigidity sequence for $S$, $\underset{i \to \infty}{\lim}\int_{\mathbb{T}} |z^{n_i}-1|^2 d \sigma_{g,S}=0$. Because $L^2$ convergence implies convergence almost everywhere along a subsequence, it follows that there exists $i_1,i_2,...$ such that $\sigma_{g,S}(\{z:\underset{j \to \infty}{\lim} z^{n_{i_j}} \to 1\})=\sigma_{g,S}(\mathbb{T})$. However, $\underset{ i \to \infty}{\lim}\int_{C}z^{n_i}\sigma_{f,T} \to 0$ 
for any measurable $C \subset \mathbb{T}$.
 This is because $\int_{C}z^{n_i}\sigma_{f,T}= \int_{\mathbb{T}}z^{n_i}\chi_C(z)\sigma_{f,T}$ and $\chi_C$ can be approximated in $L_2(\sigma_{f,T})$ by polynomials. The construction of $A$ in the previous paragraph shows that $\underset{n \in A}{\lim}\int_{\mathbb{T}}p(z)z^n d \sigma_{f,T}=0$ for any polynomial $p$. It follows that $\sigma_{g,S}$ is singular with respect to $\sigma_{f,T}$ for any $f \in H_c$ and $g \in L_0^2(\lambda)$.
\end{proof}

\begin{proof}[Proof of Theorem \ref{disjoint}] Notice that since $H_{pp}$ and $H_c$ are orthogonal and $U_T$ invariant if $g_1 \in H_{pp}$ and $g_2 \in H_c$ then $\sigma_{g_1+g_2,T}$ is $\sigma_{g_1,T}+\sigma_{g_2,T}$. It follows from Theorem \ref{hp} that any IET lying in the intersection of the full measure sets of IETs in Lemmas \ref{point spec} and \ref{cont spec} is disjoint from $T$.
\end{proof}
\begin{rem} The following observation motivates the proof. If $\mu$ and $\nu$ are probability measures on $S^1$ such that $z^{n_i} \to f$ weakly in $L_2(\mu)$ and $z^{n_i} \to g$ weakly in $L_2(\nu)$ and $f(z) \neq g(z)$ for all $z$ then $\nu$ and $\mu$ are singular.
\end{rem}
\begin{rem}\label{sequential} A possibly more checkable result follows from the above proof. Assume $A$ is a mixing sequence for $T$ (that is, $\underset{n \in A}{\lim} \, \mu(B \cap T^n(B'))=\mu(B)\mu(B')$ for all measurable $B$ and $B'$) then any $S$ having a rigidity sequence in $A$ is disjoint from $T$. Note that weak mixing transformations have mixing sequences of density 1.
\end{rem}
\begin{rem} Given a family of transformations $\mathcal{F}$ with a measure $\eta$ on $\mathcal{F}$ any $\mu$ ergodic $T \colon X \to X$ will be disjoint for $\eta$-almost every $S \in \mathcal{F}$ if:
\begin{enumerate}
\item Any sequence of density 1 is a rigidity sequence for $\eta$-almost every $S \in \mathcal{F}$.
\item $\eta(\{S \in \mathcal{F}: \alpha \text{ is an eigenvalue for } S\})=0$ for any $\alpha \neq 1$.
\end{enumerate}
Additionally, the results in the previous section show that a slightly stronger version of condition 1 and $\eta$-almost sure total ergodicity implies condition 2. Condition 1 on its own does not imply condition 2. To see this consider when $\mathcal{F}$ is the set of 1 element, rotation by $\alpha_0$.
\end{rem}

\section{Concluding remarks}
First, the proof of Corollary \ref{ue}.
\begin{proof}[Proof of Corollary \ref{ue}] This follows from Theorem \ref{disjoint}, the fact that almost every IET is uniquely ergodic (\cite{masur} and \cite{gauss}) and the following Lemma.
\end{proof}
\begin{lem} If $T$ and $S$ are uniquely ergodic with respect to $\mu$ and $\nu$ respectively then any preserved measure of $T \times S$ has projections $\mu$ and $\nu$.
\end{lem}
\begin{proof} Consider $\eta$, a preserved measure of $T \times S$. $$\eta(A \times Y)= \eta\left( (T^{-n} \times S^{-n}) (A \times Y)\right)= \eta(T^{-n}(A) \times Y).$$ Therefore, $\mu_1(A):=\eta(A \times Y)$ is preserved by $T$ and so it is $\mu$. For the other projection the proof is similar.
\end{proof}
More is true in fact, for $\bold{m}_{\mathfrak{R_1}} \times$...$\times \bold{m}_{\mathfrak{R_n}}$ almost every n-tuple of IETs $(S_1,...,S_n)$, $S_1 \times... \times S_n$ is uniquely ergodic and $S_1$ is disjoint from $S_2 \times S_3 \times...\times S_n$.

Corollary \ref{ue} has an application.
 Consider $T \times S$. In our context, unique ergodicity implies minimality, which implies uniformly bounded return time to a fixed rectangle. Therefore, if we choose a rectangle $V \subset [0,1) \times [0,1)$ then the induced map of $T \times S$ on $V$ is almost surely (in $(T,S)$ or even $S$ if $T$ is uniquely ergodic) an exchange of a finite number of rectangles. 
 To see this recall that a minimal IET $T$ is measurably isomorphic to a continuous shift dynamical system $\bar{T}$ that acts on a compact space (see \cite[Section 5]{IET}). 
 Moreover, $T=\gamma_T \circ \bar{T}$ where $\gamma_T$ is continuous and at worst a two to one map (in fact it is one to one in all but a countable number of places, the orbits of discontinuities). Unique ergodicity of $T$ implies unique ergodicity of $\bar{T}$. Likewise, if $T$ and $S$ are minimal and disjoint IETs then $\bar{T}$ and $\bar{S}$ are disjoint. It follows if they are also uniquely ergodic then $\bar{T} \times \bar{S}$ is a uniquely ergodic continuous map of a compact metric space and therefore minimal. It follows from compactness, continuity and minimality that the return time to any open set under $\bar{T} \times \bar{S}$ is bounded. The continuity of $\gamma_T \times \gamma_S$ implies that the return time to a fixed rectangle is bounded under $T \times S$.

Theorem \ref{disjoint} also strengthens Corollary \ref{eigen} because transformations are not disjoint from their factors \cite[Theorem 8.4]{glas}.
\begin{cor} No transformation is a factor of a positive measure set of IETs.
\end{cor}

A number of questions have come up in relation to the results of this paper.
\begin{ques} (Bufetov) Let $\mu$ be an ergodic measure invariant under Rauzy-Veech induction. Under what conditions is $\mu \times \mu$ almost every pair of IETs disjoint?
\end{ques}
There are atomic ergodic measures of Rauzy-Veech induction that obviously fail this. However, the fact that almost every nonrotation IET is weak mixing (proven in \cite{AF}) extends to many ergodic measures of Rauzy-Veech induction. This provides hope for extending Theorem \ref{disjoint} in these settings (see Remark \ref{sequential}). However, to replicate the arguments here one would need versions of the estimates on distortion bounds and the measure of the region that shares the same matrix of Rauzy-Veech induction.
\begin{ques} Does almost every IET with a particular permutation $\pi$ have no (or possibly only obvious) isomorphic IETs with permutation $\pi$? For instance, in the permutation $(4321)$ the IET given by length vector $(a,b,c,1-(a+b+c))$ is isomorphic to $(1-(a+b+c),c,b,a)$.
\end{ques}

Section 2 showed a particular reason for rigidity occurred fairly often for almost every IET, but could occur at any time for only a small portion of IETs. Can rigidity happen at a certain time for a larger than expected portion of IETs? The following questions occurred during conversations with Boshernitzan and Veech.

\begin{ques} Can there be a rigidity sequence for a positive measure set of IETs?
\end{ques}
\begin{ques} Can there be a particular large $n$ that is an $\epsilon$ rigidity time for a large measure set of IETs in some Rauzy class?
\end{ques}

Also, Section 2 showed that for many $R$ a set of measure at least comparable to $R^{-1}$ has an $\epsilon$ rigidity time $R$. Then next question asks if there are some times where this does not happen.
\begin{ques} Is there a sequence $R_1,R_2,...$ such that for some $\epsilon$ a set of measure at most $o(R_i^{-1})$ has an $\epsilon$ rigidity time $R_i$ ?
\end{ques}

Some outstanding questions of Veech \cite{prime} are also relevant.
\begin{ques} Is almost every IET that is not of rotation type prime? (Prime means no nontrivial measurable factors.) Does almost every IET have property S? (Property $S$ says that every ergodic self joining other than the product measure is almost everywhere one to one.) Does almost every IET that is not of rotation type have nontrivial compact subgroups in their centralizer? 
\end{ques}
\section{Acknowledgments}

I would like to thank A. Bufetov for posing this question, encouragement and many helpful conversations. Also, once I got the initial result, he encouraged me to prove more even suggesting directions to pursue. I would like to thank V.~ Bergelson, M. Boshernitzan, D. Damanik, E-H. El Abdalaoui, M. Keller, H. Kr\"uger, S. Semmes, J. P. Thouvenot and W.~ Veech for helpful conversations. I would like to thank M. Lemanczyk for pointing out relevant references. I would like to thank the referee for providing many suggestions that improved the paper. I was supported by Rice University's Vigre grant DMS-0739420 and a Tracy Thomas award while working on this paper.

\end{document}